\documentclass{article}
\usepackage{amsmath}[1996/11/01]
\usepackage{amssymb,amsthm,amsxtra}
\usepackage{epic,eepic}
\usepackage{epsfig}

\def\[#1\]{\begin{equation}#1\end{equation}}
\makeatletter
\def\beq{%
   \relax\ifmmode
      \@badmath
   \else
      \ifvmode
         \nointerlineskip
         \makebox[.6\linewidth]%
      \fi
      $$
   \fi
}
\def\eeq{%
   \relax\ifmmode
      \ifinner
         \@badmath
      \else
         $$
      \fi
   \else
      \@badmath
   \fi
   \ignorespaces
}

\def\enddisplaymath{\eeq\global\@ignoretrue}
\makeatother

\newtheorem{thm}{Theorem}
\newtheorem{cor}[thm]{Corollary}
\newtheorem{lem}[thm]{Lemma}

\theoremstyle{remark}
\newtheorem*{rem}{Remark}
\newtheorem{rems}{Remark}[thm]

\theoremstyle{definition}
\newtheorem{defn}{Definition}

\numberwithin{equation}{section}
\numberwithin{thm}{section}

\DeclareMathOperator{\Exp}{E}
\DeclareMathOperator{\Var}{Var}

\DeclareMathOperator{\fp}{fp}
\DeclareMathOperator{\ifp}{fpi}
\DeclareMathOperator{\RS}{RS}

\DeclareMathOperator{\Prob}{Pr}

\newcommand{\C}{\mathbb C}
\newcommand{\R}{\mathbb R}
\newcommand{\Z}{\mathbb Z}

\newcommand{\M}{\mathbb M}

\newcommand\psymmU{%
\begin{picture}(1,1)(0,0)%
\allinethickness{0.5pt}%
\path(0,0)(0,1)(1,1)(1,0)(0,0)%
\end{picture}}
\newcommand\psymmUU{%
\begin{picture}(1,1)(0,0)%
\allinethickness{0.5pt}%
\path(0,0)(0,1)(1,1)(1,0)(0,0)%
\put(0.5,0.5){\makebox(0,0){$\cdot$}}%
\end{picture}}
\newcommand\psymmO{%
\begin{picture}(1,1)(0,0)%
\allinethickness{0.5pt}%
\path(0,0)(0,1)(1,1)(1,0)(0,0)%
\path(0,0)(1,1)%
\end{picture}}
\newcommand\psymmS{%
\begin{picture}(1,1)(0,0)%
\allinethickness{0.5pt}%
\path(0,0)(0,1)(1,1)(1,0)(0,0)%
\path(1,0)(0,1)%
\end{picture}}
\newcommand\psymmu{%
\begin{picture}(1,1)(0,0)%
\allinethickness{0.5pt}%
\path(0,0)(0,1)(1,1)(1,0)(0,0)%
\path(0,0)(1,1)%
\path(0,1)(1,0)%
\end{picture}}

\newbox\tsymmUbox
\newbox\tsymmUUbox
\newbox\tsymmObox
\newbox\tsymmSbox
\newbox\tsymmubox
\setbox\tsymmUbox =\hbox{\kern0.75pt\setlength{\unitlength}{6pt}\psymmU \kern0.75pt}
\setbox\tsymmUUbox=\hbox{\kern0.75pt\setlength{\unitlength}{6pt}\psymmUU\kern0.75pt}
\setbox\tsymmObox =\hbox{\kern0.75pt\setlength{\unitlength}{6pt}\psymmO \kern0.75pt}
\setbox\tsymmSbox =\hbox{\kern0.75pt\setlength{\unitlength}{6pt}\psymmS \kern0.75pt}
\setbox\tsymmubox =\hbox{\kern0.75pt\setlength{\unitlength}{6pt}\psymmu \kern0.75pt}
\def\tsymmU{{\copy\tsymmUbox}}
\def\tsymmUU{{\copy\tsymmUUbox}}
\def\tsymmO{{\copy\tsymmObox}}
\def\tsymmS{{\copy\tsymmSbox}}
\def\tsymmu{{\copy\tsymmubox}}
\newbox\symmUbox
\newbox\symmUUbox
\newbox\symmObox
\newbox\symmSbox
\newbox\symmubox
\setbox\symmUbox =\hbox{\kern0.75pt\setlength{\unitlength}{4.5pt}\psymmU \kern0.75pt}
\setbox\symmUUbox=\hbox{\kern0.75pt\setlength{\unitlength}{4.5pt}\psymmUU\kern0.75pt}
\setbox\symmObox =\hbox{\kern0.75pt\setlength{\unitlength}{4.5pt}\psymmO \kern0.75pt}
\setbox\symmSbox =\hbox{\kern0.75pt\setlength{\unitlength}{4.5pt}\psymmS \kern0.75pt}
\setbox\symmubox =\hbox{\kern0.75pt\setlength{\unitlength}{4.5pt}\psymmu \kern0.75pt}
\def\symmU{{\copy\symmUbox}}
\def\symmUU{{\copy\symmUUbox}}
\def\symmO{{\copy\symmObox}}
\def\symmS{{\copy\symmSbox}}
\def\symmu{{\copy\symmubox}}

\def\tsymmg{\circledast}
\def\symmg{\circledast}

\begin{document}

\title{{\bf Symmetrized random permutations}}
\author{{\bf Jinho Baik}\footnote{
Princeton University and Institite for Advanced Study,
New Jersey, jbaik@math.princeton.edu} \ \ 
and \ \ {\bf Eric M. Rains}\footnote{AT\&T Research, New Jersey, 
rains@research.att.com}}

\date{September 23, 1999}
\maketitle



\section{Introduction}\label{sec-intro}

\begin{figure}[ht]
 \centerline{\epsfig{file=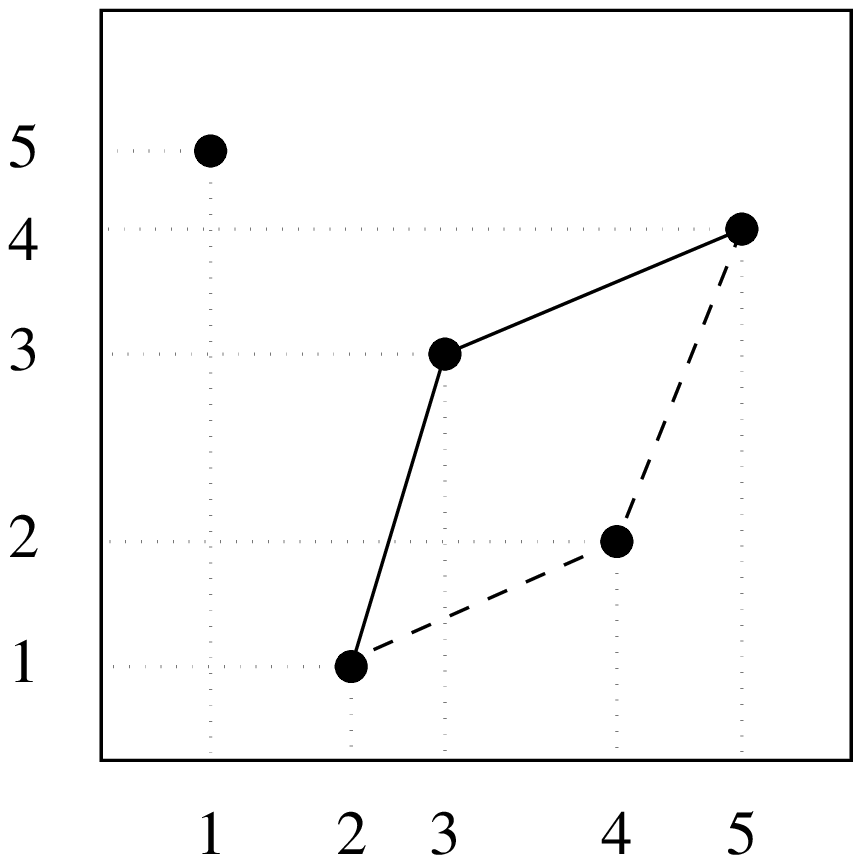, width=5cm}}
 \caption{random points in a rectangle}
\label{fig-randompoints}
\end{figure}

Suppose that we are selecting $n$ points, $p_1,p_2,\cdots,p_n$, 
at random in a rectangle, 
say $R=[0,1]\times[0,1]$ (see Figure \ref{fig-randompoints}). 
We denote by $\pi$ the configuration of $n$ random points.
With probability 1, no two points have same $x$-coordinates nor $y$-coordinates.
An up/right path of $\pi$ 
is a collection of points $p_{i_1},p_{i_2},\cdots,p_{i_k}$ 
such that $x(p_{i_1})<x(p_{i_2})<\cdots<x(p_{i_k})$ 
and $y(p_{i_1})<y(p_{i_2})<\cdots<y(p_{i_k})$.
The length of such a path is defined by the number of the points in the path.
Now we denote by $l_n(\pi)$ the length of the longest up/right path 
of a random points configuration $\pi$.

As one can see from Figure \ref{fig-randompoints}, a configuration of $n$ 
points gives rise a permutation. 
For the example at hand, the corresponding permutation is 
$\bigl( \begin{smallmatrix}1&2&3&4&5\\5&1&3&2&4\end{smallmatrix}\bigr)$.
Therefore we can identity random points in $R$ and random permutations, 
and we use the same notation $\pi$. 
In this identification, $l_n(\pi)$ is the length of the longest increasing 
subsequence of a random permutation.

The longest increasing subsequence has been of great interest 
for a long time (see e.g. \cite{AD}, \cite{OR}, \cite{BDJ}).
Especially as $n\to\infty$, it is known that 
$\Exp(l_n)\sim 2\sqrt{n}$ \cite{LS}, \cite{VK1, VK2} (also \cite{AD95,
Se95, Jo2})
and $\Var(l_n)\sim c_0n^{1/3}$ \cite{BDJ} 
with some numerical constant $c_0=0.8132\cdots$.
Moreover, the limiting distribution 
of $l_n$ after proper scaling 
is obtained in \cite{BDJ} in terms of the solution to the 
Painlev\'e II equation 
(see Section \ref{sec-result} for precise statement).
An interesting feature is that the above limiting distribution function
is also the limiting 
distribution of the (scaled) largest eigenvalue of a random GUE matrix 
\cite{TW1}, the so-called ``GUE Tracy-Widom distribution'' $F_2$. 
In other words, properly centered and scaled, the length of the longest 
increasing subsequence of a random permutation behaves statistically 
for large $n$ like the largest eigenvalue of a random GUE matrix.
There have been many papers concerning the relations on 
combinatorics and random matrix theory : 
we refer the reader to \cite{Re, Gessel, Kerov, Rains, Jo2, 
BDJ, BDJ2, TW3, Bo, 
kurtj:shape, kurtj:trans, Ok, BOO, TW:randomwords, kurtj:recent, 
BR1, BR2}. 
The purpose of this paper is to survey the analytic results of 
the recent papers \cite{BR1, BR2} 
and discuss related topics.

In random matrix theory, three ensembles play important roles,  GUE, GOE 
and GSE (see e.g. \cite{Mehta}). 
Since random permutation is related to GUE, 
it would be interesting to ask which object in combinatorics is related to 
GOE and GSE. 
For this purpose, we consider symmetrized permutations.
In terms of random points, 5 symmetry types of the rectangle $R$ 
are considered, denoted by the symbols $\tsymmU$, $\tsymmO$,
$\tsymmS$, $\tsymmUU$, and $\tsymmu$. 
Throughout this paper (and also in \cite{BR1, BR2}), the symbol 
$\tsymmg$ is used to denote an arbitrary choice of the above five possibilities.
Let $\delta=\{(t,t):0\le t\le 1\}$, the diagonal line, 
and $\delta^t=\{(t,1-t):0\le t\le 1\}$, the anti-diagonal line.
Consider the following random points selections : 
\begin{enumerate}
\item $\tsymmU$ : select $n$ points in $R$ at random.
\item $\tsymmO$ : select $n$ points in $R\setminus\delta$ 
and $m$ points in $\delta$ at random, and add 
their reflection images about $\delta$.
\item $\tsymmS$ : select $n$ points in $R\setminus\delta^t$ 
and $m$ points in $\delta^t$ at random, and add
their reflection images about $\delta^t$.
\item $\tsymmUU$ : select $n$ points at random in $R$, and add 
their rotational images about the center $(1/2,1/2)$.
\item $\tsymmu$ : select $n$ points in $R\setminus\delta$, 
$m_+$ points in $\delta$ and $m_-$ points in $\delta^t$ at random, 
and add their reflection images about both $\delta$ and $\delta^t$.
\end{enumerate}
Define the map $\iota$ on $S_n$ by $\iota(x)=n+1-x$.
Let $\fp(\pi)$ denote the number of  points satisfying $\pi(x)=x$, 
and $\ifp(\pi)$ denote the number of points satisfying $\pi(x)=\iota(x)$ 
($\ifp$ represents negated points : see Remark \ref{rem-3} below).
Each of the above process corresponds to 
picking a random permutation from each of the following ensembles :
{
\allowdisplaybreaks
\begin{align}
S^\symmU_n &= S_n\\
S^\symmO_{n,m} &= \{\pi\in S_{2n+m}:\pi=\pi^{-1},\ 
\fp(\pi)=m \}\\
S^\symmS_{n,m} &= \{\pi\in S_{2n+m}:\pi=\iota \pi^{-1}\iota ,\ 
\ifp(\pi)=m \}\\
S^\symmUU_n &= \{\pi\in S_{2n}:\pi=\iota\pi\iota\}\\
S^\symmu_{n,m_+,m_-} &= \{\pi\in S_{4n+2m_++2m_-}:
\pi=\pi^{-1},\ \pi=\iota\pi^{-1}\iota,\
\fp(\pi)=2m_+,\ \ifp(\pi)=2m_- \}.
\end{align}
}
We denote the length of the longest increasing subsequence (equivalently, 
the longest up/right path) of $\pi$ in each of the ensemble respectively by 
\begin{equation}
  L^\symmU_n,\quad L^\symmO_{n,m},\quad L^\symmS_{n,m},\quad L^\symmUU_n,\quad 
L^\symmu_{n,m_+,m_-}.
\end{equation}

\begin{rems}
  The map $\pi\mapsto \iota^{-1}\pi$ gives a bijection 
between $S^\symmO_{n,m}$ and $S^\symmS_{n,m}$. 
Thus $L^\symmS_{n,m}$ has the same statistics with the length 
of the longest \emph{decreasing} subsequence of a random 
involution with $m$ fixed points taken from $S^\symmO_{n,m}$.
From the definition, $L^\symmO_{n,m}$ is the random variable describing 
the length of the longest \emph{increasing} subsequence of a random 
involution taken from the same ensemble.
\end{rems}

\begin{rems}\label{rem-2}
  We may identify $S_{2n}$ with the set of bijections from
$\{-n,\cdots,-2,-1,1,2,\cdots,n\}$ onto itself.
In this identification, $S_n^\symmUU$ becomes the set of
signed permutations ; $\pi(x)=-\pi(-x)$.
The longest increasing subsequence problem of a random signed permutation
is considered in \cite{TW3} and \cite{Bo}.
\end{rems}

\begin{rems}\label{rem-3}
  Under the identification in Remark \ref{rem-2}, $S^\symmu_{n,m_+,m_-}$ 
becomes the set of signed involutions with $m_+$ fixed points 
and $m_-$ negated points (we call $x$ a negated point if $\pi(x)=-x$.)
\end{rems}

\bigskip
In this paper, we are interested in the statistics of $L^\symmg$ 
as $n\to\infty$.
Especially for $\tsymmO$, $\tsymmS$ and $\tsymmu$, we 
are interested in the cases when 
$m=[\sqrt{2n}\alpha]$ for $\tsymmO$, $m=[\sqrt{2n}\beta]$ for $\tsymmS$, 
and $m_+=[\sqrt{n}\alpha]$ and $m_-=[\sqrt{n}\beta]$ for $\tsymmu$ 
with fixed $\alpha, \beta \ge 0$ where $[k]$ denotes the largest 
integer less than or equal to $k$.
Then for most cases, the expected values have the same asymptotics.
Namely, if we set $N=n, 2n+m, 2n+m,2n,4n+2m_++2m_-$ for 
each of $\tsymmU, \tsymmO, \tsymmS, \tsymmUU, \tsymmu$ case respectively, 
we have 
\begin{equation}
   \lim_{N\to\infty} \frac{\Exp(L^\symmg)}{\sqrt{N}} =2, 
\end{equation}
when $0\le\alpha\le 1$ and $\beta\ge 0$ are fixed for $\tsymmO$,
$\tsymmS$ and $\tsymmu$.
When $\alpha>1$, we have different expected value in the limit 
(see Section \ref{sec-result}.)

On the other hand, the variance behaves asymptotically 
like $c_0N^{1/3}$ but now with 
different constant $c_0$ depending on the symmetry type.
It is because each symmetry type has different limiting distribution : 
$L^\symmU_n$ has GUE fluctuation, $L^\symmS_{n,m}$ GOE fluctuation 
and $L^\symmUU_n$ $\text{GUE}^2$
fluctuation 
(see Section \ref{sec-result} below for precise statements).
Here $\text{GUE}^2$ denotes the statistics of 
a superimposition of eigenvalues of two random GUE matrices.
Similarly for $\text{GOE}^2$.
The cases of $\tsymmO$ and $\tsymmu$ show more interesting features.
For $\tsymmO$, the limiting distribution function changes depending on 
the value of $\alpha=m/\sqrt{2n}$. 
The fluctuation is GSE when $\alpha<1$, GOE when $\alpha=1$ 
and Gaussian when $\alpha>1$.
By taking suitable scaling limit $\alpha\to 1$, we can find a certain smooth 
transition between GSE and GOE.
For $\tsymmu$, the value $\alpha=m_+/\sqrt{n}$ determines the limiting 
distribution ; the value $m_-$ plays no role in the transition.
The fluctuation is 
GUE when $\alpha<1$, $\text{GOE}^2$ when $\alpha=1$,
and Gaussian when $\alpha>1$.

In Section \ref{sec-dist}, we define 
the Tracy-Widom distributions for GUE, GOE and GSE as well as new classes of 
distribution functions describing the transition around $\alpha=1$.
Main results are stated in Section \ref{sec-result}, and 
Section 4 includes some applications and the related problems.
Most of the results in this article are taken from 
\cite{BR1, BR2}. 
The only new result is Theorem \ref{thm-4.2}.

\medskip
\noindent {\bf Acknowledgments.}
The authors would like to thank organizers of the workshop 
\emph{Random Matrix Models and their Applications} for their invitations.

\section{Tracy-Widom distribution functions}\label{sec-dist}

Let $u(x)$ be the solution of the Painlev\'e II (PII) equation,
\begin{equation}\label{as10}
  u_{xx}=2u^3+xu,
\end{equation}
with the boundary condition
\begin{equation}\label{as11}
  u(x)\sim -Ai(x)\quad\text{as}\quad x\to +\infty,
\end{equation}
where $Ai$ is the Airy function.
The proof of the (global) existence and the uniqueness of the solution 
was first established in \cite{HM} : 
the asymptotics as $x\to -\infty$ are 
(see e.g. \cite{HM,DZ2}) 
\begin{eqnarray}
\label{as17}  u(x) &=&
-Ai(x) + O\biggl( \frac{e^{-(4/3)x^{3/2}}}{x^{1/4}}\biggr),
\qquad\text{as $x\to +\infty$,}\\
\label{as18}  u(x) &=&
-\sqrt{\frac{-x}{2}}\biggl( 1+O\bigl(\frac1{x^2}\bigr)\biggr),
\qquad\text{as $x\to -\infty$.}
\end{eqnarray}
Recall that $Ai(x)\sim \frac{e^{-(2/3)x^{3/2}}}{2\sqrt{\pi}x^{1/4}}$
as $x\to +\infty$.
Define
\begin{equation}\label{as12}
  v(x):= \int_{\infty}^{x} (u(s))^2 ds,
\end{equation}
so that $v'(x)= (u(x))^2$.

We introduce the Tracy-Widom distributions.
(Note that $q:=-u$, which Tracy and Widom used in their papers,
solves the same
differential equation with the boundary condition
$q(x)\sim +Ai(x)$ as $x\to\infty$.)

\begin{defn}[Tracy-Widom distribution functions]\label{def1}
Set
\begin{eqnarray}
\label{as46-1}   F(x) &:=& \exp\biggl(\frac12\int_x^{\infty} v(s)ds\biggr)
= \exp\biggl(-\frac12\int_x^{\infty} (s-x)(u(s))^2ds\biggr),\\
\label{as47-1}   E(x) &:=& \exp\biggl(\frac12\int_x^{\infty} u(s)ds\biggr),
\end{eqnarray}
and set
\begin{eqnarray}
   F_2(x) :=& F(x)^2 &=\exp\biggl(-\int_x^{\infty} (s-x)(u(s))^2ds\biggr),\\
   F_1(x) :=& F(x)E(x) &= \bigl(F_2(x)\bigr)^{1/2}
e^{\frac12\int_x^{\infty} u(s)ds} ,\\
   F_4(x) :=& F(x)\bigl[E(x)^{-1}+E(x)\bigr]/2
&= \bigl(F_2(x)\bigr)^{1/2}
\biggl[ e^{-\frac12\int_x^{\infty} u(s)ds} +
e^{\frac12\int_x^{\infty} u(s)ds} \biggr] /2 .
\end{eqnarray}
\end{defn}

In \cite{TW1} and \cite{TW2},
Tracy and Widom proved that under proper centering
and scaling, the distribution of the largest eigenvalue of
a random GUE/GOE/GSE matrix converges to $F_2(x)$ / $F_1(x)$ / $F_4(x)$
as the size of the matrix becomes large.
We note that from the asymptotics \eqref{as17} and \eqref{as18},
for some positive constant $c$,
\begin{eqnarray}
\label{as46}   F(x) &=& 1+ O\bigl( e^{-cx^{3/2}}\bigr)
\qquad \text{as $x\to+\infty$,}\\
\label{as47}   E(x) &=& 1+ O\bigl( e^{-cx^{3/2}}\bigr)
\qquad \text{as $x\to+\infty$,}\\
\label{as48}   F(x) &=& O\bigl( e^{-c|x|^{3}}\bigr)
\qquad\qquad \text{as $x\to-\infty$,}\\
\label{as49}   E(x) &=& O\bigl( e^{-c|x|^{3/2}}\bigr)
\quad\qquad \text{as $x\to-\infty$.}
\end{eqnarray}
Hence in particular, $\lim_{x\to+\infty} F_\beta(x) = 1$ and
$\lim_{x\to-\infty} F_\beta(x) = 0$, $\beta=1,2,4$.
Monotonicity of $F_\beta(x)$ follows from the fact that $F_\beta(x)$
is the limit of a sequence of distribution functions.
Therefore $F_\beta(x)$ is indeed a distribution function.


\bigskip
As indicated in Introduction, we need new classes of distribution 
functions to describe the phase transitions from 
$\chi_{\text{GSE}}$ to $\chi_{\text{GOE}}$ 
and from $\chi_{\text{GUE}}$ to $\chi_{\text{GOE}^2}$.
First we consider the Riemann-Hilbert problem (RHP) for the Painlev\'e II 
equation \cite{FN, JMU}.
Let $\Gamma$ be the real line  $\R$, oriented from $+\infty$ to $-\infty$.
Let $m(\cdot\thinspace;x)$ be the solution of the following RHP :
\begin{equation}\label{as20}
 \begin{cases}
    m(z;x) \qquad \text{is analytic in $z\in\C\setminus\Gamma$,}\\
    m_+(z;x)=m_-(z;x) \begin{pmatrix} 1& -e^{-2i(\frac43z^3+xz)}\\
e^{2i(\frac43z^3+xz)}& 0 \end{pmatrix} \quad \text{for $z\in\Gamma$,}\\
    m(z;x) = I+O\bigl(\frac1{z}\bigr) \qquad \text{as $z\to\infty$.}
 \end{cases}
\end{equation}
Here $m_+(z;x)$ (resp,. $m_-$) is the limit of $m(z';x)$ as $z'\to z$
from the left (resp., right) of the contour $\Gamma$ :
$m_\pm(z;x)=\lim_{\epsilon\downarrow 0}m(z\mp i\epsilon;x)$.
Relation \eqref{as20} corresponds to 
the RHP for the PII equation with the special monodromy data
$p=-q=1, r=0$ (see \cite{FN, JMU}, also \cite{FZ,DZ2}).
In particular if the solution is expanded at $z=\infty$,
\begin{equation}\label{as7.27}
   m(z;x) = I+ \frac{m_1(x)}{z} + O\bigl(\frac1{z^2}\bigr),
\qquad \text{as $z\to\infty$},
\end{equation}
we have
\begin{eqnarray}
\label{as22}   2i(m_1(x))_{12} = -2i(m_1(x))_{21} &=& u(x), \\
\label{as23}   2i(m_1(x))_{22} = -2i(m_1(x))_{11}&=& v(x),
\end{eqnarray}
where $u(x)$ and $v(x)$ are defined in \eqref{as10}-\eqref{as12}.
Therefore the Tracy-Widom distributions above are expressed in terms of
the residue at $\infty$ of the solution to the RHP \eqref{as20}. 
It is noteworthy that the new distributions below which interpolate the 
Tracy-Widom distributions require additional information 
of the solution of RHP.

\begin{defn}\label{def3}
  Let $m(z;x)$ be the solution of RHP \eqref{as20} and
denote by $m_{jk}(z;x)$ the $(jk)$-entry of $m(z;x)$.
For $w>0$, define
\begin{equation}\label{as61}
   F^\symmO(x;w) :=
F(x) \biggl\{\bigl[m_{22}(-iw;x)-m_{12}(-iw;x) \bigr]E(x)^{-1}
+ \bigl[m_{22}(-iw;x)+m_{12}(-iw;x)\bigr]E(x) \biggr\}/2,
\end{equation}
and for $w<0$, define
\begin{equation}\label{as62}
 \begin{split}
   F^\symmO(x;w) &:=
e^{\frac83w^3-2xw}F(x) \\
&\times \biggl\{\bigl[-m_{21}(-iw;x)+m_{11}(-iw;x) \bigr]E(x)^{-1}
- \bigl[m_{21}(-iw;x)+m_{11}(-iw;x)\bigr]E(x) \biggr\}/2.
 \end{split}
\end{equation}
Also define
\begin{eqnarray}
   F^\symmu(x;w) &:=& m_{22}(-iw;x)F_2(x),
\quad\qquad\qquad\qquad w>0,\\
   F^\symmu(x;w) &:=& -e^{\frac83w^3-2xw}m_{21}(-iw;x)F_2(x),
\qquad w<0.
\end{eqnarray}
\end{defn}

\medskip
First $F^\symmO(x;w)$ and $F^\symmu(x;w)$ are real 
from Lemma \ref{lem12} (i) below.
Note that $F^\symmO(x;w)$ and $F^\symmu(x;w)$ are continuous at $w=0$
since at $z=0$, the jump condition of the RHP \eqref{as20} implies
$(m_{12})_+(0;x)= -(m_{11})_-(0;x)$ and
$(m_{22})_+(0;x)= -(m_{21})_-(0;x)$.
In fact, $F^\symmO(x;w)$ and $F^\symmu(x;w)$ are entire
in $w\in\C$ from the RHP \eqref{as20}.

From \eqref{as46}-\eqref{as49} and \eqref{as94}-\eqref{as97} below,
we see that
\begin{equation}
  \lim_{x\to+\infty} F^\symmO(x;w), F^\symmu(x;w) =1, \qquad
\lim_{x\to-\infty} F^\symmO(x;w), F^\symmu(x;w) = 0
\end{equation}
for any fixed $w\in\R$.
Also Theorem \ref{thm3} below shows 
that $F^\symmO(x;w)$ and $F^\symmu(x;w)$
are limits of distribution functions, implying that
they are monotone in $x$.
Therefore, $F^\symmO(x;w)$ and $F^\symmu(x;w)$ are indeed
distribution functions for each $w\in\R$.


We close this section summarizing some properties of $m(-iw;x)$ in the 
following lemma. 
In particular the lemma implies
that $F^\symmO(x;w)$ interpolates between $F_4(x)$ and $F_1(x)$, 
and $F^\symmu(x;w)$ interpolates between $F_2(x)$ and $F_1(x)^2$ 
(see Corollary \ref{cor16}).

\begin{lem}\label{lem12}
  Let $\sigma_3=\bigl(\begin{smallmatrix} 1&0\\0&-1 \end{smallmatrix}
\bigr)$, $\sigma_1=\bigl(\begin{smallmatrix} 0&1\\1&0 \end{smallmatrix}
\bigr)$, and set $[a,b]=ab-ba$.
\begin{enumerate}
\item 
  For real $w$, $m(-iw;x)$ is real.
\item
For fixed $w\in\R$, we have
\begin{eqnarray}
\label{as94}  m(-iw;x) &=& \bigl( I+\bigl( e^{-cx^{3/2}} \bigr) \bigr)
\begin{pmatrix}
1 & -e^{\frac83w^3-2xw} \\ 0 & 1
\end{pmatrix}, \qquad  \text{$w>0$, $x\to+\infty$,}\\
\label{as95}  m(-iw;x) &=& \bigl( I+\bigl( e^{-cx^{3/2}} \bigr) \bigr)
\begin{pmatrix}
1 & 0 \\ -e^{-\frac83w^3+2xw} & 1
\end{pmatrix}, \qquad  \text{$w<0$, $x\to+\infty$,}\\
\label{as96}  m(-iw;x) &\sim& \frac1{\sqrt{2}} \biggl( \begin{smallmatrix}
1&-1\\ 1&1 \end{smallmatrix} \biggr)
e^{(-\frac43w^3+xw)\sigma_3}
e^{(\frac{\sqrt{2}}3(-x)^{3/2}+\sqrt{2}w^2(-x)^{1/2})\sigma_3},
\quad  \text{$w>0$, $x\to-\infty$,}\\
\label{as97}  m(-iw;x) &\sim& \frac1{\sqrt{2}} \biggl( \begin{smallmatrix}
1&1\\ -1&1 \end{smallmatrix} \biggr)
e^{(-\frac43w^3+xw)\sigma_3}
e^{(-\frac{\sqrt{2}}3(-x)^{3/2}-\sqrt{2}w^2(-x)^{1/2})\sigma_3},
\quad  \text{$w<0$, $x\to-\infty$.}
\end{eqnarray}
\item 
  For any $x$, we have
\begin{equation}\label{as98}
  \lim_{w\to 0^+} m(-iw;x) =
\lim_{w\to 0^-} \sigma_1 m(-iw;x)\sigma_1  = \begin{pmatrix}
\frac12\bigl(E(x)^2+E(x)^{-2}\bigr) & -E(x)^2 \\
\frac12\bigl(-E(x)^2+E(x)^{-2}\bigr) & E(x)^2
\end{pmatrix}.
\end{equation}
\item 
   For fixed $w\in\R\setminus\{0\}$,
$m(-iw;x)$ solves the differential equation
\begin{equation}\label{as62.5}
   \frac{d}{dx}m = w[m,\sigma_3]+ u(x)\sigma_1 m,
\end{equation}
where $u(x)$ is the solution of the PII equation \eqref{as10}, \eqref{as11}.
\end{enumerate}
\end{lem}

\begin{cor}\label{cor16}
   We have
\begin{eqnarray}
   F^\symmO(x;0) &=& F_1(x),\\
   \lim_{w\to\infty} F^\symmO(x;w) &=& F_4(x), \\
   \lim_{w\to-\infty} F^\symmO(x;w) &=& 0, \\
\label{as116}   F^\symmu(x;0) &=& F_1(x)^2, \\
   \lim_{w\to\infty} F^\symmu(x;w) &=& F_2(x), \\
   \lim_{w\to-\infty} F^\symmu(x;w) &=& 0.
\end{eqnarray}
\end{cor}

\begin{proof}
   The values at $w=0$ follow from \eqref{as98}.
For $w\to\pm\infty$, note that from the RHP \eqref{as20},
we have $\lim_{z\to\infty} m(z;x)= I$.
\end{proof}

\section{Main Results}\label{sec-result}

As in the Introduction, let $N$ denote 
$n, 2n+m, 2n+m,2n,4n+2m_++2m_-$ for
each of $\tsymmU, \tsymmO, \tsymmS, \tsymmUU, \tsymmu$ case respectively.
We scale the random variables : for permutations and involutions, 
\begin{equation}
\chi^\symmU_n = \frac{L^\symmU_n-2\sqrt{N}}{N^{1/6}}, \qquad 
\chi^\symmO_{n,m} = \frac{L^{\symmO}_{n,m}-2\sqrt{N}}{N^{1/6}}, \qquad 
\chi^\symmS_{n,m} = \frac{L^{\symmS}_{n,m}-2\sqrt{N}}{N^{1/6}}, 
\end{equation}
and for signed permutations and signed involutions, 
\begin{equation}
\chi^\symmUU_{n} = \frac{L^\symmUU_n-2\sqrt{N}}{2^{2/3}N^{1/6}}, \qquad 
\chi^\symmu_{n,m_+,m_-} =
\frac{L^{\symmu}_{n,m_+,m_-}-2\sqrt{N}}{2^{2/3}N^{1/6}}.
\end{equation}
All the results in this section are taken from \cite{BR2} which utilizes 
the algebraic work of \cite{BR1}.

First, we state the results for random permutations and 
random signed permutations.
Random permutations show GUE fluctuation in the limit, 
while random signed permutations have $\text{GUE}^2$ fluctuation.
\begin{thm}\label{thm1}
  For fixed $x\in\R$, 
\begin{align}
  \lim_{n\to\infty} \Prob\bigl( \chi^\symmU_n \le x\bigr) 
&= F_2(x), \\
  \lim_{n\to\infty} \Prob\bigl( \chi^\symmUU_n \le x\bigr) 
&= F_2(x)^2.
\end{align}
\end{thm}

For the involution cases, we have the following limits.
\begin{thm}\label{thm2}
 For each fixed $\alpha$ and $\beta$, and for fixed $x\in\R$, 
we have : for $\tsymmO$, 
\begin{align}
  \lim_{n\to\infty} \Prob\bigl( \chi^\symmO_{n,[\sqrt{2n}\alpha]}
\le x \bigr) &= F_4(x), 
\qquad 0\le \alpha<1, \\
  \lim_{n\to\infty} \Prob\bigl( \chi^\symmO_{n,[\sqrt{2n}]}
\le x \bigr) &= F_1(x), \\
  \lim_{n\to\infty} \Prob\bigl( \chi^\symmO_{n,[\sqrt{2n}\alpha]}
\le x \bigr) &= 0, 
\qquad\qquad \alpha>1, 
\end{align}
for $\tsymmS$, 
\begin{align}
  \lim_{n\to\infty} \Prob\bigl( \chi^\symmS_{n,[\sqrt{2n}\beta]}
\le x \bigr) &= F_1(x), 
\qquad \beta\ge 0, 
\end{align}
and for $\tsymmu$, 
\begin{align}
  \lim_{n\to\infty} 
\Prob\bigl( \chi^\symmu_{n,[\sqrt{n}\alpha],[\sqrt{n}\beta]}
\le x \bigr) &= F_2(x), 
\qquad 0\le\alpha<1, \beta\ge 0, \\
  \lim_{n\to\infty} 
\Prob\bigl( \chi^\symmu_{n,[\sqrt{n}],[\sqrt{n}\beta]}
\le x \bigr) &= F_1(x)^2, 
\qquad \beta\ge 0, \\
  \lim_{n\to\infty} 
\Prob\bigl( \chi^\symmu_{n,[\sqrt{n}\alpha],[\sqrt{n}\beta]}
\le x \bigr) &= 0, 
\qquad\qquad \alpha>1, \beta\ge 0.
\end{align}
\end{thm}

This theorem shows that for $\tsymmO$ and $\tsymmu$, 
the limiting distributions differ depending on the value of $\alpha$.
As indicated earlier in the Introduction, as $\alpha\to 1$ at a certain rate, 
we obtain smooth transitions. 
From Corollary \ref{cor16},
the following results are consistent with Theorem \ref{thm2}.
\begin{thm}\label{thm3}
For fixed $w\in\R$, $\beta\ge 0$ and $x\in\R$, 
\begin{align}
  \lim_{n\to\infty} \Prob\bigl( \chi^\symmO_{n,m}
\le x \bigr) &= F^\symmO(x;w),
\qquad m=[\sqrt{2n} - 2w(2n)^{1/3}],\\
     \lim_{n\to\infty}
\Prob\bigl( \chi^\symmu_{n,m_+,m_-}
\le x \bigr) &= F^\symmu(x;w),
\quad\qquad m_+=[\sqrt{n} - 2wn^{1/3}],\ m_-=[\sqrt{n}\beta].
\end{align}
\end{thm} 

When $\alpha>1$, Theorem \ref{thm2} shows that we have used 
inappropriate scaling. 
In a proper scaling, we obtain normal distribution $N(0,1)$.
\begin{thm}\label{thm4}
  For fixed $\alpha>1$ and $\beta\ge 0$, as $n\to\infty$, 
\begin{align}
  \frac{L^\symmO_{n,[\sqrt{2n}\alpha]}-(\alpha+1/\alpha)\sqrt{N}}
{\sqrt{(1/\alpha-1/\alpha^3)}N^{1/4}} 
&\to N(0,1) \qquad \text{in distribution},\\
  \frac{L^\symmu_{n,[\sqrt{n}\alpha],[\sqrt{n}\beta]}
-(\alpha+1/\alpha)\sqrt{N}}
{\sqrt{(1/\alpha-1/\alpha^{3})}N^{1/4}} 
&\to N(0,1) \qquad \text{in distribution}.
\end{align}
\end{thm}

All the above results are on the convergence in distribution.
We also have convergence of moments for all the cases.

\begin{thm}\label{thm5}
  For each case of the above theorems, all the moments of 
the random variable converge to the moments of the corresponding 
limiting distribution.
\end{thm}

From this result, we can obtain the asymptotics of variances.
Especially for $\tsymmU$, the variance is 
\begin{equation}
  \lim_{N\to\infty} \frac{\Var(l^\symmU_n)}{N^{1/3}}
=\int_{-\infty}^\infty x^2 dF_2(x) -
\biggl( \int_{-\infty}^\infty xdF_2(x)\biggr)^2.
\end{equation}
Evaluating the integrals, we obtain the value $0.8132\cdots$ 
(see \cite{TW1}).

The outline of the proofs is as follows.
First we consider the Poisson generating function. 
It is to let the number of points be Poisson.
For example, we define 
\begin{equation}
   P^\symmU_l(t)=e^{-t^2} \sum_{n=0}^{\infty} \frac{t^{2n}}{n!}
\Prob\bigl( L^\symmU_n\le l \bigr).
\end{equation}
The de-Poissonization lemma \cite{Jo2} tells us that 
in the limit $n\to\infty$, $\Prob\bigl( L^\symmU_n\le l \bigr) 
\sim P^\symmU_l(n^2)$, hence it is enough to obtain the asymptotics of 
the Poisson generating function.
The crucial point is that there is a determinantal formula for 
each Poisson generating function.
For the cases $\tsymmU$ and $\tsymmUU$, the determinant is 
of Toeplitz type \cite{Gessel, Rains}, while for the rest, 
it is of Hankel type \cite{BR1}.
In fact, as in \cite{Gessel}, there are general identities between 
sum of Schur functions and determinantal formulae (see \cite{BR1} 
for details), which can be used to consider other type of 
Young tableaux problems (see Subsection \ref{subsec-Jo} below).
Now general theory connects 
Toeplitz/Hankel determinants and orthogonal polynomials.
It turns out that to analyze all the above cases, only one set of 
orthogonal polynomials are needed, namely 
the orthogonal polynomials on the unit circle with respect to the 
weight $e^{2t\cos\theta} d\theta/(2\pi)$.

Now following Fokas, Its and Kitaev \cite{FIK}, there is a Riemann-Hilbert
representation for orthogonal polynomials.
Let $\Sigma$ be the unit circle in the complex plane oriented
counterclockwise. 
Let $Y(z)$ be a $2\times 2$ matrix-valued function satisfying 
\begin{equation}
\begin{cases}
   Y(z) \quad\text{is analytic in $\C\setminus\Sigma$},\\
   Y_+(z)=Y_-(z) \begin{pmatrix} 1&\frac{1}{z^k}e^{t(z+z^{-1})}\\
0&1 \end{pmatrix}, \quad z\in\Sigma, \\
   Y(z) \bigl(\begin{smallmatrix} z^{-k}&0\\0&z^k
\end{smallmatrix} \bigr) = I+O(z^{-1})  
\qquad \text{as $z\to\infty$}, 
\end{cases}
\end{equation}
where $Y_+(z)$ (resp. $Y_-(z)$) denotes the limit of $Y(z')$ 
as $z'\to z$ satisfying $|z'|<1$ (resp. $|z'|>1$).
Then one finds that for example, the 11 entry of $Y(z)$ is the 
$k^{\text{th}}$ monic orthogonal polynomial with respect to the weight 
$e^{2t\cos\theta} d\theta/(2\pi)$.
Once we have a Riemann-Hilbert representation, we can employ the
steepest-descent method (Deift-Zhou method) developed by Deift and Zhou
\cite{DZ1} to find asymptotics as parameters become large (or small).
For our case, the parameters are $t$ and $k$, and taking proper scaling,
we obtain precise asymptotics which eventually yield the convergence in 
distribution and convergence of moments.
We also mention that in our analysis, equilibrium measures play a
crucial role as in the papers \cite{DKMVZ, DKMVZ2, DKMVZ3}.

\section{Applications and related topics}\label{sec-appl}

\subsection{Random involutions and random signed involutions}

The ensemble $S^\symmO_{n,m}$ is the set of involutions with 
$n$ 2-cycles and $m$ 1-cycles.
In the previous section, we considered the limiting statistics when 
$n$ and $m$ are related by $m=[\sqrt{2n}\alpha]$ with $\alpha$ 
being finite ; either fixed or $\alpha\to 1$ with certain rate.
It is of interest to consider the whole set of involutions 
without constraints on the number of fixed points.
Similarly, the signed involutions without constraint 
on the number of fixed points and negated points is also of interest.
We define the ensembles of involutions and signed involutions 
\begin{align}
  \tilde{S}_n &=\{\pi\in S_n : \pi=\pi^{-1} \}, \\
  \tilde{S}_n^\symmu &= \{\pi\in S_{2n} : 
\pi=\pi^{-1},\ \pi=\iota\pi\iota \}, 
\end{align}
and denote by $\tilde{L}_n(\pi)$ and $\tilde{L}^\symmu_n(\pi)$ the length
of the longest increasing subsequence of $\pi\in\tilde{S}_n$ 
and that of $\pi\in\tilde{S}_n^\symmu$, 
respectively.

Noting $\tilde{S}_n=\bigcup_{2k+m=n} S^\symmO_{k,m}$, 
we have 
\begin{equation}\label{e-appl-1}
  \Prob\bigl(\tilde{L}_n\le l\bigr)= \frac{1}{|\tilde{S}_{k,m}|} 
\sum_{2k+m=n} \Prob\bigl(L^\symmO_{k,m}\le l\bigr) |S_{k,m}|. 
\end{equation}
It is not difficult to check that (see pp.66-67 of \cite{Kn}) 
as $n\to\infty$, the main contribution to the sum 
$|\tilde{S}_{k,m}|=\sum_{2k+m=n} |S_{k,m}|$ comes from 
$\sqrt{2k}-(2k)^{\epsilon+1/4}\le m\le\sqrt{2k}+(2k)^{\epsilon+1/4}$.
Comparing with the scaling $m=[\sqrt{2k}-2w(2k)^{1/3}]$ 
in Theorem \ref{thm3}, the main contribution to the sum 
\eqref{e-appl-1} comes from when $w=0$, or $\alpha=1$.
Thus we obtain GOE Tracy-Widom distribution function in the limit.
Similarly, we can obtain the convergence of moments.
The signed involution case is analogous.

\begin{thm}
For fixed $x\in\R$, 
\begin{align}
  \lim_{n\to\infty} 
\Prob\biggl( \frac{\tilde{L}_n-2\sqrt{n}}{n^{1/6}} \le x\biggr)
&= F_1(x), \\
\lim_{n\to\infty} 
\Prob\biggl( \frac{\tilde{L}_n^\symmu-2\sqrt{n}}{2^{2/3}n^{1/6}}\le x\biggr)
&= F_1(x)^2.
\end{align}
We also have convergence of all the moments.
\end{thm}

This result should be compared with the results on random permutation 
and random signed permutation where the limiting distribution 
was $F_2(x)$ and $F_2(x)^2$ under the same scaling of the above 
(see \cite{BDJ}, \cite{TW3}, \cite{Bo} and Theorem \ref{thm1} above).

\begin{rem}
  As in the permutation case, the length of the longest \emph{increasing} 
subsequence and the length of the longest \emph{decreasing} subsequence 
of random involutions have the same statistics. 
This can be seen by noting that there is a bijection 
(called Robinson-Schensted correspondence, see e.g. \cite{Kn}) 
between 
the set $\tilde{S}_n$ of involutions of $n$ letters 
and the set of standard Young tableaux 
of size $n$, and the rows and the columns of standard Young tableaux 
have the same statistics under the push forward of the uniform probability 
distribution on $\tilde{S}_n$ under this bijection.
\end{rem}

\subsection{$\beta$-Plancherel measure on the set of Young diagrams}

Let $Y_n$ be the set of Young diagrams, or equivalently partitions, 
of size $n$.
Given a partition $\lambda=(\lambda_1,\lambda_2,\cdots) \vdash n$, 
let $d_\lambda$ denote the number of standard Young tableaux of shape
$\lambda$. 
We introduce the $\beta$-Plancherel measure $M^\beta_n$ on $Y_n$ 
defined by 
\begin{equation}
   M^\beta_n(\lambda)=\frac{d_\lambda^\beta}
{\sum_{\mu\vdash n}d_\mu^\beta}, \qquad \lambda\in Y_n.
\end{equation}
When $\beta=2$, this is the Plancherel measure which arises in the
representation theory.
We are interest in the typical shape and the fluctuation of $\lambda$ 
where $\lambda$ is taken randomly from the probability space $Y_n$  
with $M^\beta_n$.

A motivation introducing the above measure is the result of 
Regev \cite{Re}. 
In \cite{Re}, it is proved that for fixed $\beta >0$ and fixed $l$,
as $n\to\infty$,
\begin{equation}\label{int-as-1}
   \sum_{\substack{\lambda\vdash n \\ \lambda_1\le l}}
\bigl(d_\lambda\bigr)^\beta
\sim \biggl[ \frac{l^{l^2/2}l^n}{(\sqrt{2\pi})^{(l-1)/2}n^{(l-1)(l+2)/4}}
\biggr]^\beta \frac{n^{(l-1)/2}}{l!}
\int_{\R^l} e^{-\frac12\beta l\sum_{j}x_j^2}
\prod_{j<k} |x_j-x_k|^\beta d^l x.
\end{equation}
The multiple integral on the right hand side is 
the Selberg integral which can be computed exactly for each $\beta$.
When $\beta=1,2,4$, this integral is 
the normalization constant of the probability 
density of eigenvalues in 
GOE, GUE and GSE, respectively (see e.g. \cite{Mehta}).
So the basic question is if the $\beta$ in the definition 
of the $\beta$-Plancherel measure corresponds to 
the $\beta$ in the random matrix theory.

The well-known Robinson-Schensted correspondence \cite{Sc} 
establishes a bijection between 
$S_n$ and the pairs of standard Young tableaux with the same
shape of size $n$, $\RS : \pi \mapsto (P(\pi),Q(\pi))$.
Especially, we obtain $\sum_{\mu\vdash n} d_\mu^2= |S_n|=n!$.
Moreover, under $\RS$, the length of the longest increasing subsequence 
of $\pi\in S_n$ is equal to the number of boxes in the first 
row of $P(\pi)$ (or equally of $Q(\pi)$).
Therefore under $\RS$, the Plancherel measure $M^2_n$ is 
simply the push forward of 
the uniform probability measure on $S_n$ to $Y_n$, 
and the number of boxes in the first row of a random Young diagram 
and the length of the longest increasing subsequence of a random 
permutation have the same statistics : GUE ($\beta=2$) 
fluctuation in the limit.

If $\RS(\pi)=(P,Q)$, then $\RS(\pi^{-1})=(Q,P)$ (see e.g. \cite{Kn}).
Therefore the set of involutions $\tilde{S}_n$ is bijective to 
the set of (single) standard Young tableaux, and 
the number of boxes in the first row of a random Young diagram 
taken under the probability $M^1_n$ has the same statistics 
with the length of the longest increasing subsequence of a random
involution : GOE ($\beta=1$) fluctuation in the limit.

In fact, it is shown in \cite{Ok, BOO, kurtj:recent} that, 
for the case of $\beta=2$, in the large $n$ limit, 
the number of boxes in the $k^{\text{th}}$ of a random Young diagram 
corresponds to the $k^{\text{th}}$ largest eigenvalue 
of a random GUE matrix.
Also the typical shape of $\lambda$ is obtained in \cite{LS} and 
\cite{VK1} which is related to Wigner's semicircle law \cite{Kerov}.
On the other hand, the second row of a random Young diagram 
for $\beta=1$ is 
discussed in \cite{BR2} implying that it corresponds to 
the second eigenvalue of a random GOE matrix.
It would be interesting to obtain similar results for general row 
for $\beta=1$ and also for general $\beta$.

\subsection{Random turn vicious walker model}

Random permutations and random involutions arise also in certain random 
walk process in a one dimensional integer lattice.
We call a particle left-movable (resp. right-movable) 
if its left (resp. right) site is vacant.
Initially, $p$ particles are located at the points $1,2,3,\cdots,p$.
At each time step $t$ ($t=1,2,\cdots$), 
one particle 
among the left-movable particles 
is selected at random,  
and is moved to its left site
(so at $t=1$, the leftmost particle is moved.)
Suppose this process is repeated for $n$ time steps. 
It is found that \cite{forrester99}
there is a bijection between all the possible configurations 
and the set of Young tableaux with at most $p$ rows of size $n$ : 
simply, in the $k^{\text{th}}$ row, write the times 
when the particle is originally at the point $k$ (see Figure \ref{fig-walk}).
\begin{figure}[ht]
 \centerline{\epsfig{file=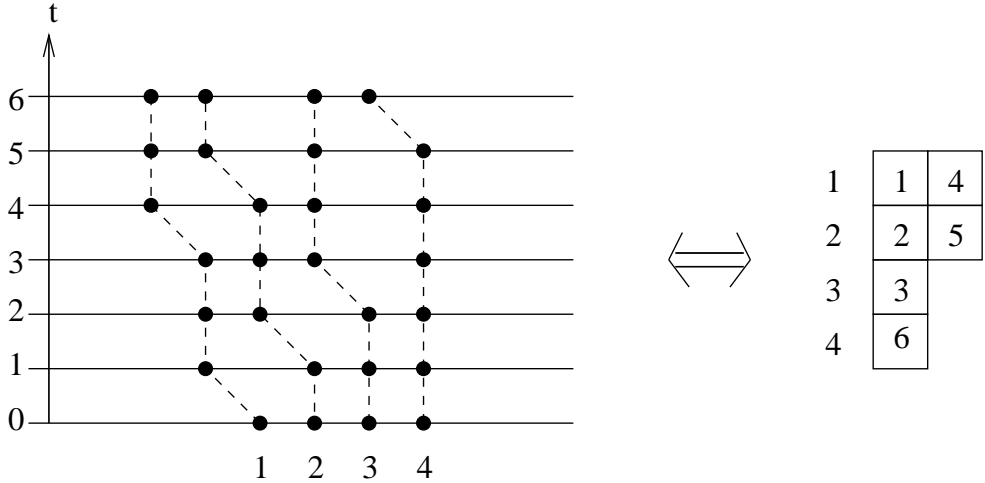, width=13cm}}
 \caption{random turn vicious walkers and Young tableaux}
\label{fig-walk}
\end{figure}
In this correspondence the number of moves made by the $k^{\text{th}}$ 
particle (counted from the left) is the number of boxes 
in the $k^{\text{th}}$ row of the corresponding Young tableau.
If we take the limit $p\to\infty$, we remove the constraint 
on the number of rows. 
In other words, we start with countably many particles located at 
$1,2,3,\cdots$ and at each 
time step, we move a randomly selected left-movable particle 
to its left site.
Then the statistics of the number of moves made by the leftmost particle 
in $n$ time steps 
is identical to that of the number of boxes 
in the first row of a random Young 
diagram under the probability $M^1_n$ (or the length of the longest 
increasing subsequence of a random involution).
Hence in the large $n$ limit, we obtain GOE fluctuation.

Now suppose after $n$ steps of left moves, 
by taking $n$ steps of right moves, 
we want the particles brought back to their original positions.
The first $n$ steps gives a Young tableaux and the next $n$ step gives another 
Young tableaux.
So we are in the situation of pairs of Young tableaux, i.e. $\beta=2$ case.
Especially if the number of particles are infinite, 
the statistics of the moves made by the leftmost 
particle is identical to the length of the longest increasing 
subsequence of a random permutation : GUE fluctuation
\cite{forrester99}.

\subsection{Symmetrized versions of Johansson's model}\label{subsec-Jo}

In \cite{kurtj:shape}, Johansson introduced a certain model 
which has several probabilistic interpretations, 
as a randomly growing Young diagram, a totally asymmetric one dimensional 
exclusion process, a certain zero-temperature directed polymer in a random 
environment or as a kind of first-passage site percolation model.
Here we consider the symmetrized versions.
The limiting distributions are parallel to the results of Section 
\ref{sec-result} depending on the symmetry type.

Each model is defined on $\M$, a subset of $\Z^2$, and at each site 
$(i,j)\in\M$, we define a random variable $w(i,j)$ of geometric
distribution. 
We are interested in ``the length of the longest up/right path''.
Let $g(q)$ denote the geometric distribution with parameter $q$.
Let $0<q<1$, and let $\alpha,\beta\ge 0$ such that 
$\alpha \sqrt{q}, \beta \sqrt{q}<1$.
Otherwise stated, the random variables $w(i,j)$ are independent each other.
We denote by $(j,k)\nearrow (j',k')$ the set of up/right paths $\pi$ from 
$(j,k)$ to $(j',k')$.

\begin{enumerate}
\item Let $\M=\Z_+^2$. 
Let 
\begin{equation}
   w(i,j)\sim g(q).
\end{equation}
Define 
\begin{equation}
  G^\symmU(N)=\max \{\sum_{(i,j)\in\pi} w(i,j) : \pi\in (1,1)\nearrow(N,N) \}.
\end{equation}
\item Let $\M=\Z^2$. 
Let $w(i,j)=0$ if $i=0$ or $j=0$, and for $i,j\neq 0$, let 
\begin{equation}
   w(i,j)=w(-i,-j)\sim g(q).
\end{equation}
Define 
\begin{equation}
   G^\symmUU(N)=\max \{ \sum_{(i,j)\in\pi} w(i,j) : 
\pi\in (-N,-N)\nearrow(N,N) \}.
\end{equation}
\item Let $\M=\Z_+^2$. 
Let 
\begin{eqnarray}
   w(i,j)=w(j,i) &\sim& g(q),\quad   i\neq j, \\
   w(i,i)&\sim& g(\alpha\sqrt{q}).
\end{eqnarray}
Define 
\begin{equation}
  G^\symmO(N)=\max \{ \sum_{(i,j)\in\pi} w(i,j) : \pi\in (1,1)\nearrow(N,N) \}.
\end{equation}
\item Let $\M=\Z_+\times \Z_-$.
Let 
\begin{eqnarray}
  w(i,-j)=w(j,-i) &\sim& g(q), \quad i\neq -j, \\
  w(i,-i) &\sim& g(\beta\sqrt{q}).
\end{eqnarray}
Define 
\begin{equation}
   G^\symmS(N)=\max \{ \sum_{(i,j)\in\pi} w(i,j) : 
\pi\in (1,-N)\nearrow (N,-1) \}.
\end{equation}
\item Let $\M=\Z^2$.
Let $w(i,j)=0$ if $i=0$ or $j=0$. Otherwise 
\begin{eqnarray}
   w(i,j)=w(-i,-j) &\sim& g(q),\quad |i|\neq|j|, \\
   w(i,i)=w(-i,-i) &\sim& g(\alpha\sqrt{q}), \\
   w(i,-i)=w(-i,i) &\sim& g(\beta\sqrt{q}).
\end{eqnarray}
Define 
\begin{equation}
   G^\symmu(N)=\max \{ \sum_{(i,j)\in\pi} w(i,j) : 
\pi\in(-N,-N)\nearrow (N,N) \}.
\end{equation}
\end{enumerate}

We are interested in the limiting distributions of $G^\symmg(N)$.
\begin{thm}\label{thm-4.2}
Set 
\begin{equation}
  \eta(q)=\frac{2\sqrt{q}}{1-\sqrt{q}}, 
\qquad \rho(q)=\frac{q^{1/6}(1+\sqrt{q})^{1/3}}{1-\sqrt{q}}.
\end{equation}
We have 
\begin{enumerate}
\item 
\begin{equation}\label{e-appl-2}
  \lim_{N\to\infty}
\Prob\biggl(\frac{G^\symmU(N;q)-\eta(q)N}
{\rho(q)N^{1/3}}\le x\biggr)
= F_2(x).
\end{equation}
\item 
\begin{equation}
  \lim_{N\to\infty}
\Prob\biggl(\frac{G^\symmUU(N;q)-\eta(q)(2N)}
{2^{2/3}\rho(q)(2N)^{1/3}}\le x\biggr)
= F_2(x)^2.
\end{equation}
\item 
\begin{eqnarray}
  \lim_{N\to\infty}
\Prob\biggl(\frac{G^\symmO(N;q,\alpha)-\eta(q)N}
{\rho(q)N^{1/3}}\le x\biggr)
&=& F_4(x), \qquad 0\le\alpha<1, \\
\label{e-appl-3}
  \lim_{N\to\infty}
\Prob\biggl(\frac{G^\symmO(N;q,\alpha)-\eta(q)N}
{\rho(q)N^{1/3}}\le x\biggr)
&=& F^\symmO(x;w), \quad \alpha=1-\frac{2w}{\rho(q)N^{1/3}}, \\
  \lim_{N\to\infty}
\Prob\biggl(\frac{G^\symmO(N;q,\alpha)-\eta(q)N}
{\rho(q)N^{1/3}}\le x\biggr)
&=& 0, \qquad\qquad \alpha>1.
\end{eqnarray}
\item
\begin{equation}
  \lim_{N\to\infty}
\Prob\biggl(\frac{G^\symmS(N;q,\beta)-\eta(q)N}
{\rho(q)N^{1/3}}\le x\biggr)
= F_1(x), \qquad 0\le\beta.
\end{equation}
\item 
\begin{eqnarray}
  \lim_{N\to\infty}
\Prob\biggl(\frac{G^\symmu(N;q,\alpha,\beta)-\eta(q)(2N)}
{2^{2/3}\rho(q)(2N)^{1/3}}\le x\biggr)
&=& F_2(x), \qquad 0\le\alpha<1, \beta\ge 0, \\
  \lim_{N\to\infty}
\Prob\biggl(\frac{G^\symmu(N;q,\alpha,\beta)-\eta(q)(2N)}
{2^{2/3}\rho(q)(2N)^{1/3}}\le x\biggr)
&=& F^\symmu(x;w), \quad \alpha=1-\frac{2w}{\rho(q)(2N)^{1/3}}, \beta\ge 0 \\
  \lim_{N\to\infty}
\Prob\biggl(\frac{G^\symmu(N;q,\alpha,\beta)-\eta(q)(2N)}
{2^{2/3}\rho(q)(2N)^{1/3}}\le x\biggr)
&=& 0, \qquad\qquad \alpha>1, \beta\ge 0.
\end{eqnarray}
\end{enumerate}
\end{thm}

\begin{rem}
  The results \eqref{e-appl-2} and \eqref{e-appl-3} with $w=0$ 
are obtained in \cite{kurtj:recent}. 
For $\tsymmU$, 
the longest up/right path in $(1,1)\nearrow (M,N)$, $M\neq N$, 
is also considered.
\end{rem}

The proof of the above theorem is analogous to that of 
the result in Section \ref{sec-result}.
There again is a Toeplitz/Hankel determinantal formula for each 
case. 
In fact, the only change with the symmetrized permutation case is 
that we have $(1+q+2\sqrt{q}\cos\theta)^N d\theta/(2\pi)$ instead of 
$e^{2t\cos\theta} d\theta/(2\pi)$ for the weight. 
In analyzing the asymptotics of orthogonal polynomials, we need
different scaling, but once we scale, the analysis is parallel.
One important common property of the above two weights is 
that both of the supports of their equilibrium measures change 
from the full circle to a part of the circle with one gap 
depending on $2t/k<1$ and $2t/k>1$, and $N/k<\eta(q)^{-1}$ 
and $N/k>\eta(q)^{-1}$. 
And also in the one gap case, 
the equilibrium measures decay like a square root at the ends 
of the supports.

\bibliographystyle{plain}
\bibliography{ref}

\end{document}